\documentstyle[12pt, twoside]{article}
\input amssym.def
\input amssym.tex

\newenvironment{demo}[1]%
{\vskip-\lastskip\medskip
  \noindent
  {\em #1.}\enspace
  }%
{\qed\par\medskip
  }

\newcommand{\qed}{
  \strut\hfill
  \mbox{$\Box$}
  }

{\vskip-\lastskip
  \begin{trivlist}    \item[]
      {\bf Axiom #1 }
      }%
{\end{trivlist}
  }

\newtheorem{theorem}{Theorem}[section]

\newtheorem{lemma}{Lemma}[section]

\newtheorem{example}{Example}[section]

\newtheorem{remark}{Remark}[section]

\newtheorem{proposition}{Proposition}[section]

\newcommand{\B}{\cal B}

\newcommand{\C}{ \Bbb C}

\newcommand{\End}{ \mbox{End} }

\newcommand{\fm}{ {}^m {\cal F} }

\newcommand{\fn}{ {}^n {\cal F} }

\newcommand{\fnd}{ {}^n {\cal F}_{\bf d} }

\newcommand{\Hom}{ \mbox{Hom }}

\newcommand{\hf}{ \frac12}

\newcommand{\Mm}{ {}^m {\cal M} }

\newcommand{\Mn}{ {}^n {\cal M} }

\newcommand{\Nn}{ {\cal N}_n }

\newcommand{\Nnm}{  {\cal N}_{\min (n,m)} }

\newcommand{\Nk}{  {\cal N}_k }

\newcommand{\pair}{(gl_n, gl_m) }

\newcommand{\Vm}{ {}^m {V}_{\lambda} }

\newcommand{\Vn}{ {}^n {V}_{\lambda} }

\newcommand{\Vk}{ {}^k {V}_{\lambda} }

\newcommand{\Wm}{ {}^m {W} }

\newcommand{\Wn}{ {}^n {W} }

\newcommand{\znn}{ {}^n {Z}^n }

\newcommand{\zdd}{ {}^d {Z}^d }

\newcommand{\zkm}{ {}^k {Z}^m }

\newcommand{\znk}{ {}^n {Z}^k }

\newcommand{\znd}{ {}^n {Z}^d }

\newcommand{\znm}{ {}^n {Z}^m }

\newcommand{\zmm}{ {}^m {Z}^m }

\newcommand{\zmk}{ {}^m {Z}^k }

\newcommand{\Z}{\Bbb Z}

\begin{document}
\title{
Lagrangian construction of the $\pair$-duality
  }
\author{
  Weiqiang Wang
}

\date{}
\maketitle

\begin{abstract}
 We give a geometric realization of
 the symmetric algebra of the tensor space $\C^n \bigotimes \C^m$
 together with the action of the dual pair $\pair$
 in terms of lagrangian cycles in the cotangent bundles of
 certain varieties. We establish geometrically
 the equivalence  between the $\pair$-duality and Schur duality.
 We establish the connection between
 Springer's construction of (representations of) Weyl groups
 and Ginzburg's construction
 of (representations of) Lie algebras of type $A$.
\end{abstract}

\setcounter{section}{-1}
\section{Introduction}

  The left action of $gl_n$ (resp. $gl_m$)
on the first (resp. second) factor of the
tensor product $\C^n \bigotimes \C^m$ induces natural left
actions on the $d$-th symmetric tensor space
$S^d (\C^n \bigotimes \C^m)$. These two actions clearly
commute with each other and they form a dual pair in the
sense of Howe (cf. \cite{H}). Following Howe, we have the
isotypic decomposition
\begin{eqnarray}  \label{eq_howe}
  S^d (\C^n \bigotimes \C^m)
  =  \bigoplus_{\lambda \in {\cal P}^d_{\min (n,m)} }
   \Vn \bigotimes \Vm
\end{eqnarray}
where ${\cal P}^d_k$ is the set of partitions of $d$ into at most $k$ parts
and $\Vn$ (resp. $\Vm$) is the irreducible
module of $gl_n$ (resp. $gl_m$) with highest
weight $\lambda$, and $\min (n,m)$ denotes the minimum of $n$ and $m$.

On the other hand, we consider the set $\fn$ of $n$-step flags
${\frak F} = ( 0 = F_0 \subset \ldots \subset F_n = \C^d )$ of the
vector space $\C^d$ of complex dimension $d$, with the induced
action of the general linear group $GL_d$. In such a setup,
Beilinson, Lusztig and MacPherson \cite{BLM} constructed the
quantum group for $gl_n$. Inspired by their construction, Ginzburg
\cite{G} obtained a micro-local version of it for the enveloping
algebra $U(gl_n)$ in terms of lagrangian cycles in the cotangent
bundles of $\fn \times \fn$. This is an analog of the Springer
theory for Weyl groups (cf. e.g. \cite{CG, Hu}). Note that
statements are made in \cite{G, CG} in terms of $sl_n$ although
the construction indeed yields $gl_n$. For our purpose, it is
important to stick to $gl_n$.

Let $\cal N$ be the nilpotent cone in the general linear Lie
algebra $gl_d$ and let $\Nn$ be the subset of
$n$-step nilpotents in $\cal N$. Let $\Mn$ be the set
$$
 \Mn :=
  \{(x, {\frak F}) \in \Nn \times \fn \mid x (F_i ) \subset F_{i -1},
     i = 1, \dots, n \}.
$$
One of the key varieties which
we introduce in this paper is the fibred product
$$\znm := \Mn \times_{\Nnm} \Mm \subset \Mn \times \Mm.
$$
This generalizes the early considerations (cf. \cite{BLM, G})
of the variety $\znn$ in our notation.
As is shown in \cite{G, CG}, $H(\znn)$ admits a canonical
associative algebra structure and there exists a surjective algebra
homomorphism $\rho_n : U(gl_n) \longrightarrow H(\znn)$,
where $H(Z)$ denotes the (component-wise) top
dimensional Borel-Moore homology of the variety $Z$.

Our first main result is to give a geometric realization
of the $\pair$-duality (\ref{eq_howe}) in terms of the
variety $\znm$,
which can be summarized by the following commuting
diagram (Theorem~\ref{th_lagrdual}):
\begin{eqnarray*}
 \begin{array}{ccccc}
     U(gl_n) & \looparrowright
             & S^d (\C^n \bigotimes \C^m) & \looparrowleft & U(gl_m) \\
  \downarrow &
             & \parallel                  &  & \downarrow \\
     H(\znn) & \looparrowright
             & H(\znm)                    & \looparrowleft & H(\zmm) \\
  \parallel  &
             & \parallel                  &  & \parallel   \\
  \bigoplus_{\lambda \in {\cal P}^d_n }\End (\Vn) & \looparrowright
             & \bigoplus_{\lambda \in {\cal P}^d_{\min (n,m)} }
                 \Hom ( \Vm, \Vn )
             & \looparrowleft &  \bigoplus_{\lambda \in {\cal P}^d_m}\End (\Vm)
 \end{array}
\end{eqnarray*}
Here $\looparrowright$ and $\looparrowleft$ denote left and
right algebra actions.

One advantage of our new approach to the $\pair$-duality is that a
natural basis analogous to the canonical basis in quantum groups
is given by the fundamental classes of the irreducible components
of $\znm$. We can easily identify each isotypic component inside
$H(\znm)$. It will be of great interest to construct other Howe's
duality (cf. \cite{H}) in the spirit of this paper. Such a
construction will shed some light on the geometric construction of
($q$-deformed) enveloping algebras of other classical Lie
algebras.

The Springer theory for the symmetric group $S_d$, regarded as the
Weyl group of $gl_d$, is built on the setup of the $GL_d$-diagonal
action on ${\B} \times {\B}$, where $\B$ denotes the variety of
all complete flags in $\C^d$. Grojnowski and Lusztig \cite{GL}
considered the fibred product $\Wn := \fn \times_{\Nn} \B$ and
provided a way to obtain the q-deformed Schur duality on the space
of  $d$-th tensors of $\C^n$ by studying $\Wn$ (also see
\cite{GRV} for further applications).

We formulate the classical Schur duality
in terms of Borel-Moore homology. Our second main
result is to show how
the interplay among cycles in $\Wn, Z$ and $\znm$ gives rise to
the equivalence between the $\pair$-duality and Schur duality.
We shall also derive Springer's theorem on the Weyl groups
of type $A $ from the $\pair$-duality. Along the way,
we clarify the precise relations between Springer's
original construction of representations of Weyl groups
and Ginzburg's construction of representations of
Lie algebras of type $A$.

The results of this paper can be generalized to the quantum group
setting based on constructions over finite fields,
or in terms of equivariant K-theory (cf. \cite{BLM, GL, GRV, CG}).

The plan of this paper is as follows. In Sect.~\ref{sec_basic}
we review the convolution in Borel-Moore homology and
introduce some new varieties we need later.
In Sect.~\ref{sec_main} we present the lagrangian construction
of the $\pair$-duality. In Sect.~\ref{sec_schur} we establish
the equivalence between the $\pair$-duality and Schur duality.
\section{The variety $\znm_k$}
   \label{sec_basic}
   Given a locally compact space $Z$ the Borel-Moore
homology $H_{\bullet} (Z)$ with complex coefficient
is defined to be the ordinary (relative) homology of the
pair $(\widehat{Z}, \infty )$ where $\widehat{Z}$ is the
one-point comactification of $Z$. See \cite{CG}
for more on the Borel-Moore homology. We denote by
$H(Z)$ the subspace of $H_{\bullet} (Z)$ consisting of
(component-wise) top homology classes. In all of our applications,
each connected component of $Z$ is of pure dimension, and
so $H(Z)$ is spanned by the
fundamental classes of all irreducible components of $Z$.

Given two closed subsets $Z$ and $Z^{'}$ of a smooth oriented
manifold $M$ of real dimension $s$, we can define an intersection pairing
$$
  \cap : H_i (Z) \times H_j (Z^{'}) \longrightarrow
    H_{i +j -s} (Z\cap Z^{'} ).
$$

Now let $M_1, M_2$ and $M_3$ be smooth oriented manifolds
of real dimension $m_1, m_2$ and $m_3$ respectively. Let
$Z_{12} \subset M_1 \times M_2$, $Z_{23} \subset M_2 \times M_3$
be closed subsets. Assume that the map
$$
 p_{13} : p_{12}^{-1} (Z_{12}) \cap p_{23}^{-1} (Z_{23})
    \longrightarrow M_1 \times M_3
$$
to be proper. Define the set-theoretic composition
\begin{eqnarray*}
  Z_{12} \circ Z_{23} =
  \{ (x_1, x_3) \in M_1 \times M_3 \mid \mbox{there exists }
    x_2   \nonumber  \\
  \quad\quad\quad\quad
  \mbox{ such that } (x_1, x_2) \in Z_{12} \mbox{ and }
    (x_2, x_3) \in Z_{23}
  \} .
\end{eqnarray*}

We observe that $Z_{12} \circ Z_{23}$ is the image of
$p_{12}^{-1} (Z_{12}) \cap p_{23}^{-1} (Z_{23})$ under the
projection $p_{13}$. We define a {\em convolution}
\begin{eqnarray*}
  \star : H_i (Z_{12}) \times H_j (Z_{23})
                   \longrightarrow H_{i +j -m_2} (Z_{12} \circ Z_{23})
\end{eqnarray*}
by letting
$$
  (c_{12}, c_{23})
   \mapsto  c_{12} \star c_{23} = (p_{13})_*
     (c_{12} \boxtimes [M_3] \cap [M_1] \boxtimes c_{23}).
$$
$M_1, M_3$ are allowed to be disconnected. Here
$[M_1]$ and $ [M_3]$ are understood as the sum of
the fundamental classes of connected components of
$M_1$ and $M_3$.

It is easy to check that the
convolution in the Borel-Moore homology satisfies the
associativity. Furthermore it has an important
{\em dimension property}: let
$$
  p = \frac{m_1 + m_2}2, q = \frac{m_2 + m_3}2,
  r = \frac{m_1 + m_3}2.
$$
Then the convolution induces a map
$ H_p (Z_{12}) \times H_q (Z_{23}) \longrightarrow
  H_r (Z_{12} \circ Z_{23}).
$

Now let us introduce various geometric objects we will
need in the subsequent sections.
We fix a positive integer $d$ throughout this paper.
We denote by $\fn$ the set of all $n$-step partial flags
$$
  \frak F = (0 = F_0 \subset F_1 \subset \ldots \subset F_n = \C^d ).
$$ $\fn$ is a smooth compact manifold with connected components
$\fnd$ parameterized by the set of all maps ${\bf d}: [1, n]
\longrightarrow \Z_+$ such that $\sum_{i =1}^n d_i = d$, where
$d_i$ denotes the value of ${\bf d}$ at $i$. Here $$
 \fnd =
 \{\frak F = (0 = F_0 \subset F_1 \subset \ldots \subset F_n = \C^d )
 \mid \dim F_i / F_{i -1} = d_i \}.
$$

Denote by $\cal N$ the nilpotent cone in $gl_d$ and by
$\Nn$ the subset of $n$-step nilpotents in $\cal N$, namely
$$
 \Nn =
 \{ x \in \End (\C^d) \mid x^n = 0 \}.
$$
Clearly $\Nn$ is a closed subvariety of $\cal N$. We define
$$
 \Mn :=
  \{(x, \frak F) \in \Nn \times \fn \mid
   x (F_i) \subset F_{i -1}, i = 1, \ldots, n \}
$$
and denote its connected component by
$\Mn_{\bf d} := \Mn \cap (\Nn \times \fn_{\bf d})$.
$GL_d$ acts on $\Mn, \Nn$ and $\fn$. This makes the following natural maps
induced from the projections $GL_d$-equivariant:
\begin{eqnarray*}
 \begin{array}{rcl}
    & \Mn &   \\
  \quad \mu\swarrow& & \searrow\pi \quad  \\
  \Nn \;\; & & \fn
 \end{array}
\end{eqnarray*}

One has a natural isomorphism of vector bundles between $\Mn$ and the
cotangent bundle $T^* (\fn)$ which makes the following
diagram commute:
\begin{eqnarray*}
 \begin{array}{rcl}
  \Mn \quad  & \approx &  T^* (\fn) \\
   \;\;\pi\searrow& & \swarrow \widetilde{\pi} \\
  & \fn &
 \end{array}
\end{eqnarray*}
Here $\widetilde{\pi}$ is the natural projection from $T^* (\fn )$
to the base manifold $\fn$. We shall denote ${\fn}_x : =
\mu^{-1}(x)$ where $\mu :\Mn \rightarrow \Nn$ and $x\in \Nn$.

For integers $n, m >0$ and $k \geq 0$, we introduce the following
variety which will play a fundamental role in our geometric
construction of $\pair$-duality:
\begin{eqnarray*}
 \znm_k & :=& \Mn \times_{\Nk} \Mm \\
      & =& \{ (x, \frak F, {\frak F}^{'})
          \in \Nk \times \fn \times \fm  \\
      & & \quad\quad
       \mid x (F_i ) \subset F_{i -1},
            x (F_j^{'} ) \subset F_{j -1}^{'},
      1 \leq i \leq n, 1 \leq j \leq m \}.
\end{eqnarray*}

\begin{remark} \rm
  Clearly $\znm_k$ is a closed subvariety of $ \znm_{k +1}$, and
 $\znm_k$ stabilizes and is equal to
 the variety $\znm : = \Mn \times_{\cal N} \Mm$ when
 $k \geq \min (n,m)$. We will often implicitly use the identities
 $ \znm_{\min (n, k, m)} =  \znm_{\min (n, k)}
  =  \znm_{\min (k, m)} = \znm_k$.
\end{remark}

We note that
\begin{eqnarray}  \label{eq_cotang}
   \Mn \times \Mm
       \approx T^* (\fn) \times T^* (\fm)
           \stackrel{sign}{\approx}  T^* (\fn \times \fm).
\end{eqnarray}
The second isomorphism above is given by changing the sign
of the symplectic structure in the second factor $T^* (\fm)$ as usual.

$\znm$ is not connected in general. However we have the following
crucial dimension property component-wise.
Let ${\znm_k}({\cal O})$ be the preimage of a nilpotent
 $GL_d$-orbit $\cal O$ under the projection from $\znm_k$ to $\Nk$.

\begin{proposition}   \label{prop_mid}
  Let ${\bf d}_1$ (resp. ${\bf d}_2$) be a
partition of $d$ in ${\cal P}^d_n$ (resp. ${\cal P}^d_m$).
Let $Z^{\alpha}$ be an irreducible component of $\znm_k$
 contained in $\Mn_{{\bf d}_1} \times \Mm_{{\bf d}_2}$. Then we have
 $$
   \dim Z^{\alpha} = \hf \dim
        \left(  \Mn_{{\bf d}_1} \times \Mm_{{\bf d}_2}
        \right).
 $$
\end{proposition}

\begin{demo}{Proof}
 Clearly $\znm_k$ is the union of ${\znm_k}({\cal O})$ over
 all nilpotent conjugacy classes in $\Nk$. We can write
 $$
   {\znm_k}({\cal O}) = GL_d \times_{G(x)} (\fn_x \times \fm_x ),
 $$
 where $x \in \cal O$ and $ G(x) $ is the stabilizer of $x$ in $GL_d$.
 Thus an irreducible component $Z_{\alpha}$ of $\znm_k$ is the closure
 of a unique component of ${\znm_k}({\cal O})$
 (say in the connected component $\fn_{ {\bf d}_1} \times \fm_{ {\bf d}_2}$)
 of the form
 $$
    GL_d \times_{G(x)}
   (\fn_x^{\alpha}
    \times \fm_x^{\beta}),
 $$
 where $\fn^{\alpha}_x$ (resp. $\fm^{\beta}_x$)
 is an irreducible component in $\fn_x \cap \fn_{ {\bf d}_1}$
 (resp. $\fm_x \cap \fm_{ {\bf d}_2}$).
 According to a theorem of Spaltenstein (cf. \cite{S}),
 the variety $\fn_x \cap \fn_{ {\bf d}_1}$
 is  connected and of pure dimension so that
 \begin{eqnarray*}
  \dim {\cal O}_x + 2 \dim (\fn_x \cap \fn_{ {\bf d}_1})
   = 2 \dim \fn_{{\bf d}_1}.
 \end{eqnarray*}
 Similarly we have
 \begin{eqnarray*}
 \dim {\cal O}_x + 2 \dim (\fm_x \cap \fm_{ {\bf d}_2})
    = 2 \dim \fm_{ {\bf d}_2}.
 \end{eqnarray*}
 It follows by comparing with Eq.~(\ref{eq_cotang}) that
 \begin{eqnarray*}
 \dim Z^{\alpha} & =& \dim {\znm_k}({\cal O}) \\
  & =& \dim {\cal O}_x + \dim ( \fn_x \cap \fn^{\alpha}_{ {\bf d}_1})
   + \dim ( \fm_x \cap \fm^{\beta}_{ {\bf d}_2})  \\
  & =&  \dim \fn_{{\bf d}_1} + \dim \fm_{ {\bf d}_2}  \\
  & =&  \hf \dim  \left(  \Mn_{{\bf d}_1} \times \Mm_{{\bf d}_2}
        \right).
 \end{eqnarray*}
 This completes the proof.
\end{demo}

In the case when $k \geq \min (n,m)$ and so $\znm_k = \znm$,
we have the following description of $\znm$ which
strengthens Proposition~\ref{prop_mid}. The proof
is the same as for a similar statement in Springer theory
(cf. Proposition~3.3.4, \cite{CG}).

\begin{proposition}   \label{prop_lagrang}
 The variety $\znm$ is a union of the conormal bundles to all
 the $GL_d$-orbits in $\fn \times \fm$.
 Each irreducible component is the closure of the conormal
 bundle to a unique $GL_d$-orbit.
\end{proposition}
\begin{remark}  \rm  \label{rem_lag}
  It follows that $\znm$ is a lagrangian subvariety of $\Mn \times \Mm$. So
 is $\znm_k$ since $\znm_k$ is a subvariety of $\znm$ and
 $\znm_k$ has half the dimension of $\Mn \times \Mm$
 component-wise by Proposition~\ref{prop_mid}. Note that
 in general $\znm_k$ is not a union of some conormal bundles to all
 the $GL_d$-orbits in $\fn \times \fm$.
 $\znm_0 = \fn \times \fm$ for $k =0$ is such an example.
\end{remark}
\section{Lagrangian construction of $\pair$-duality}
     \label{sec_main}
  We make the following observation. The case when $n =m$ is well known
(cf. e.g. \cite{BLM, CG}).

\begin{lemma}  \label{lem_corresp}
  The orbits of the diagonal action of $GL_d$ on
 $\fn \times \fm$ are parameterized by the $n \times m$ matrices with
 non-negative integral entries and with the sum of all entries equal to $d$.
\end{lemma}

The correspondence is defined as follows:
given a pair of flags $({\frak F}, {\frak F^{'} }) \in \fn \times \fm$,
we let the $(i,j)$-th entry of the $n \times m$ matrices to be
the dimension of the quotient
 $        F_i \cap F_j^{'} /
          (F_{i -1} \cap F_j^{'} +
            F_i \cap F_{j -1}^{'}) ,
$
$1 \leq i \leq n, 1 \leq j \leq m$. In particular, it
follows from the lemma that the number of $GL_d$-diagonal orbits
in $\fn \times \fm$ is equal to $\dim S^d (\C^n \bigotimes \C^m )$.

 One easily verifies the set-theoretic composition
\begin{eqnarray*}
  \znm_a \circ \zmk_b = \znk_{\min (a,m,b) },
  \quad a, b \geq 0, \;n, m, k >0.
\end{eqnarray*}
By applying the convolution in Borel-Moore homology introduced
in Sect.~\ref{sec_basic} to our situation, we obtain a map
\begin{eqnarray*}
  H_{\bullet}(\znm_a ) \times H_{\bullet}(\zmk_b)
    \rightarrow H_{\bullet}(\znk_{\min (a,m,b) } ).
\end{eqnarray*}
Thanks to Proposition~\ref{prop_mid} and the dimension
property of the Borel-Moore homology, we obtain
the following proposition.

\begin{proposition}  \label{prop_gener}
 The convolution induces a map
 \begin{eqnarray}  \label{eq_convol}
   \star :
   H (\znk_a ) \times H (\zkm_b) \rightarrow H (\znm_{\min (a,k,b) }).
 \end{eqnarray}
\end{proposition}

We list below some important consequences of Proposition~\ref{prop_gener}.
We recall that $\znm_r = \znm$ if $ r \geq \min (n,m)$.

\begin{proposition}    \label{prop_spec}
 \begin{enumerate}
  \item[1)] $(H (\znn_r), \star)$  is an associative algebra with unit.
 In particular, $(H (\znn), \star)$  is an associative algebra.
  \item[2)] The algebra $H (\znn)$ acts on $H (\znm_k )$ from the left
 while $H (\zmm)$ acts on $H (\znm_k)$ from the right by convolution.
 These two actions commute with each other.
 \end{enumerate}
\end{proposition}

\begin{demo}{Proof}
  Putting $a =b=r$ and $ k =m =n$
 in (\ref{eq_convol}) we obtain the associative algebra structure on
 $H (\znn_r)$. The unit is given by the fundamental class of
 the whole $\znn_r$. Putting $a = k =n$ and replacing $b$ by $k$
 in (\ref{eq_convol}) gives us a left
 $H(\znn)$-action on $H(\znm_k)$. Similarly we obtain a right
 $H(\zmm)$-action on $H(\znm_k)$ by setting $b =k =m$ and
 replacing $a$ by $k$. It is clear by construction these two actions
 on $H (\znm)$ commute with each other.
\end{demo}

\begin{remark}  \rm
 The statement that $(H (\znn), \star)$  is an associative algebra
 is due to Ginzburg \cite{G}.
\end{remark}

We denote by ${\cal P}^d$ the set of partitions
of $d$ and by ${\cal P}^d_n$ the set of partitions
of $d$ into at most $n$ parts.
Note that ${\cal P}_{n}^d$ is in one-to-one
correspondence with the set of nilpotent conjugacy classes in
the $n$-step nilpotent cone $\Nn $:
given a partition $\lambda \in {\cal P}^d_n$, let
${\lambda}^t = (a_1, a_2, \ldots )$
be the transpose of $\lambda$
$(n \geq a_1 \geq a_2 \geq \ldots \geq 0)$; we associate to $\lambda$ the
Jordan form $x_{\lambda}$ in $\Nn$ consisting of Jordan blocks of
size $a_1, a_2, \ldots$.
We denote by $c_{\lambda}$ the nilpotent conjugacy classes of $x_{\lambda}$.

One easily shows that
$$
 \znn \circ \fn_x = \fn_x, \quad \fm_x \circ \zmm = \fm_x.
$$
Following Ginzburg \cite{G, CG},
$H (\fn_x)$ is an irreducible module over $H(\znn)$.
Given an algebra $A$ and a left $A$-module $V$, we endow
$V^{\vee} := \Hom (V, \C)$ a right $A$-module structure by
letting $(\check{v} \cdot a) (v) = \check{v} (a \cdot v)$, where
$\check{v} \in V^{\vee}, a \in A, v \in V$. It can be shown
(cf. \cite{CG}) that $H({\fm}_x)_R$ is isomorphic to
$ H(\fm_x)^{\vee}$ as a right $H(\znn)$-module, here
$R$ denotes right module.

\begin{theorem}   \label{th_duality}
   Under the commuting left action of $H(\znn)$ and the right
 action of $ H(\zmm)$, the space
 $H(\znm_k)$ decomposes as follows:
 $$
  H (\znm_k) = \bigoplus_{\lambda \in {\cal P}^d_{\min (n,k,m)} }
             Hom ( H( {\fm}_{x_{\lambda}}), H( {\fn}_{x_{\lambda}}) ).
 $$
 Here $Hom ( H( {\fm}_{x_{\lambda}}), H( {\fn}_{x_{\lambda}}) )$
 is understood as $H( {\fn}_{x_{\lambda}}) \bigotimes
 H( {\fm}_{x_{\lambda}})^{\vee}$.
\end{theorem}

Note that in the case when $k =m = n$, we recover a theorem of
Ginzburg (see Theorem 4.1.23 in \cite{CG}).
Our theorem in the general case can be proved by following the same
line as in \cite{CG} which we sketch below.

\begin{demo}{Proof}
  We introduce a partial order on the set of nilponent orbits
 ${\cal O} \subset \cal N$:
 $$
  {\cal O}^{'} \leq {\cal O} \Longleftrightarrow
  {\cal O}^{'} \subset \overline{\cal O} , \quad
  {\cal O}^{'} < {\cal O} \Longleftrightarrow
  {\cal O}^{'} \subset \overline{\cal O} ,\; {\cal O}^{'} \neq {\cal O},
 $$
 where $\overline{\cal O}$ is the closure of $\cal O$. Put
 $$
 A_{\leq \cal O} = \bigsqcup_{ {\cal O}^{'} \leq {\cal O} }
   \znm_k ( {\cal O}^{'}  ), \quad
 A_{< \cal O} = \bigsqcup_{ {\cal O}^{'} < {\cal O} }
   \znm_k ( {\cal O}^{'}  ).
 $$
 These are closed subvarieties of $\znm_k$. We observe that
 $ H( A_{< \cal O} ) \subset H(A_{\leq \cal O} )$
 as a module of $ H(\znn )$ and $ H(\zmm) )$.
 Denote $H_{\cal O} = H(A_{\leq \cal O} ) / H( A_{< \cal O} ).$
 Then $H_{\cal O}$ is spanned by the fundamental classes of
 irreducible components of $\znm_k ( {\cal O})$.
 One can show that
 as a $ (H(\znn ), H(\zmm) )$-module there is an isomorphism
 between $H_{\cal O}$ and $H( \fn_x ) \bigotimes H( \fm_x)_R$.

 The partial order $\leq$ gives a filtration of $H(\znm_k )$ by the
 submodule $H( A_{\leq \cal O})$. The graded space associated to
 this filtration as
 a $ (H(\znn ), H(\zmm) )$-module is isomorphic to $H(\znm_k)$
 since $H(\znn )$ and $ H(\zmm)$
 are semisimple algebras \cite{CG}. Therefore
 \begin{eqnarray*}
  H(\znm_k) & \approx & \bigoplus_{{\cal O} \subset \cal N} H_{\cal O} \\
    & \approx &
 \bigoplus_{\lambda \in {\cal P}^d_{\min (n,k,m)} }
             H( {\fn}_{x_{\lambda}}) \bigotimes H( {\fm}_{x_{\lambda}})_R \\
    & \approx &
 \bigoplus_{\lambda \in {\cal P}^d_{\min (n,k,m)} }
             H( {\fn}_{x_{\lambda}}) \bigotimes
             H( {\fm}_{x_{\lambda}})^{\vee} \\
    & \approx &
 \bigoplus_{\lambda \in {\cal P}^d_{\min (n,k,m)} }
             Hom ( H( {\fm}_{x_{\lambda}}), H( {\fn}_{x_{\lambda}}) ).
 \end{eqnarray*}
\end{demo}

The following theorem relates $H(\znn)$ to the enveloping
algebra $U (gl_n)$.

\begin{theorem} \label{th_hwt} {\em (Ginzburg)}
  There exists a canonical surjective
 algebra homomorphism $\rho_n : U(gl_n) \longrightarrow H(\znn)$.
 The $H(\znn)$-module $H( {\cal F}_{x_{\lambda}}) $,
 regarded as a left $gl_n$-module, has highest weight $\lambda$.
\end{theorem}

It is known that there exists an anti-involution $\omega$ on $H(\zmm)$
induced from sending $(x, {\frak F}, {\frak F}^{'})$ in $\zmm$
to $(x, {\frak F}^{'}, {\frak F})$ which is compatible via $\rho_m$ with
the Cartan anti-involution on $gl_m$ by
taking the transpose (cf. \cite{CG}). By combining the right
action of $H(\zmm)$ and the anti-involution $\omega$,
we obtain a left $H(\zmm)$-action on $H(\znm_k)$.

Given a left (resp. right) $gl_m$-module $V$, we define a right (resp. left)
$gl_m$ action on $V$ by letting an element act by its transpose.
Denote by $V^t$ the  right (resp. left) $gl_m$-module thus obtained.
As a consequence of the existence of a non-degenerate invariant form on $\Vm$
the right $gl_m$-module ${\Vm}^t$ is isomorphic to $\Vm^{\vee}$.
Thus we can reformulate Theorem~\ref{th_duality} as follows.

\begin{theorem}  \label{th_gldual}
     Under the left  action of $gl_n$ and the right action
 of $gl_m$, the space $H(\znm)$ decomposes as follows:
 $$
  H (\znm_k) = \bigoplus_{\lambda \in {\cal P}^d_{\min (n,k,m)} }
   \Vn \bigotimes \Vm^{\vee}.
 $$
\end{theorem}

\begin{remark} \rm  \label{rem_weight}
  Denote by ${\frak h}_n$ the Cartan subalgebra of diagonal
 matrices in $gl_n$ and ${}^n e_i$ $(i =1, \ldots, l)$
 the standard basis of ${\frak h}_n$ with $1$ in the $(i,i)$-th entry and
 $0$ elsewhere. By using the explicit formula
 for the homomorphism $\rho_n$ (resp. $\rho_m$) (cf. \cite{CG}),
 we can easily show that any fundamental class of an irreducible
 component of $\znm$ is a weight vector with respect to
 the action of ${\frak h}_n$ and ${\frak h}_m$. More precisely,
 given an irreducible component $Z_{A}$ of $\znm$ corresponding
 to an $n \times m$ matrix $A = (a_{ij})$ in Lemma~\ref{lem_corresp}.
 Define $\lambda_i = \sum_{j =1}^m a_{ij},$ $\mu_j = \sum_{i =1}^n a_{ij}$.
 Then $(\lambda_1, \ldots, \lambda_n)$ is the weight for
 ${\frak h}_n \subset gl_n$ (with respect to the standard
 basis $ ({}^n e_i)_{i =1}^n$) and
 $ (\mu_1, \ldots, \mu_m)$ is the weight for ${\frak h}_m \subset gl_m$.
\end{remark}

\begin{example}  \rm
 When $m =1$, the variety ${}^1 {\cal F}$ is a single point
 and so ${}^n Z^1_k = \fn$. Note that $\fn$ is also the fiber
 $\mu^{ -1} (0)$. Theorem~\ref{th_gldual} shows that as a $gl_n$-module
 $H(\fn )$ is irreducible and isomorphic to $S^d (\C^n)$.
 On the other hand, if $d =1$ then $\fn$ is a discrete set of
 $n $ points and so $\znm$ is a set of $nm$ points.
 Therefore  $H( \znm ) \approx \C^n \bigotimes \C^m$.
\end{example}

The natural left actions of $gl_n$ and $gl_m$ on
$\C^n \bigotimes \C^m$ induces left actions on the
$d$-th symmetric tensor $S^d (\C^n \bigotimes \C^m)$.
The decomposition (\ref{eq_howe}) of the space $S^d (\C^n \bigotimes \C^m)$ is
exactly the same as $H (\znm)$ in Theorem~\ref{th_gldual}
under the action of the dual pair $\pair$.
Thus we have obtained a left $\pair$-module isomorphism
\begin{eqnarray}   \label{eq_isom}
  H (\znm) \approx S^d (\C^n \bigotimes \C^m).
\end{eqnarray}
The natural isomorphism above seems to be new even in the
case when $n =m$.

It follows from Theorem~\ref{th_gldual} that
we can identify the algebra $H(\znn)$ with a direct sum
of the endomorphism algebra $\sum_{\lambda \in {\cal P}^d_n } \End ( \Vn)$.
The left (resp. right) action of $H(\znn)$
(resp. $H(\zmm)$) on $H (\znm)$ corresponds to
the natural left (resp. right) action of
$\End ( \Vn)$ (resp. $\End ( \Vm)$) on $\Vn \bigotimes \Vm^{\vee}$.
For the convenience of the reader we summarize the main
results of this section as follows.
 The right action of $gl_m$ on $S^d (\C^n \bigotimes \C^m)$ in the diagram
 is related to the left action of $gl_m$ in Eq.~(\ref{eq_howe})
 by the Cartan anti-involution in $gl_m$.

\begin{theorem} \label{th_lagrdual} {\em ($\pair$-duality)}
  We have the following commutative diagram:
\begin{eqnarray*}
 \begin{array}{ccccc}
     U(gl_n) & \looparrowright
             & S^d (\C^n \bigotimes \C^m) & \looparrowleft & U(gl_m) \\
  \downarrow &
             & \parallel                  &  & \downarrow \\
     H(\znn) & \looparrowright
             & H(\znm)                    & \looparrowleft & H(\zmm) \\
  \parallel  &
             & \parallel                  &  & \parallel   \\
  \bigoplus_{\lambda \in {\cal P}^d_n }\End (\Vn) & \looparrowright
             & \bigoplus_{\lambda \in {\cal P}^d_{\min (n,m)} }
                 \Hom ( \Vm, \Vn )
             & \looparrowleft &  \bigoplus_{\lambda \in {\cal P}^d_m}\End (\Vm)
 \end{array}
\end{eqnarray*}
 Here $\looparrowright$ and
 $\looparrowleft$ denote left and right algebra actions.
\end{theorem}

Following the same line of the proof of Theorem~\ref{th_gldual},
we can have the following identification (which is
$\pair$-equivariant):
 \begin{eqnarray*}
  \begin{array}{ccc}
   H (\znk_a )  & =&
     \bigoplus_{\lambda \in {\cal P}^d_{\min (n,a, k)} }
             \Hom(\Vk, \Vn) \\
     \bigoplus  & &\bigoplus  \\
     H (\zkm_b)&  = &
  \bigoplus_{\lambda \in {\cal P}^d_{\min (k,b,m )} }
             \Hom(\Vm, \Vk)                    \\
   \downarrow &  & \downarrow \\
  H (\znm_{\min (a,k,b ) })
           & = &
  \bigoplus_{\lambda \in {\cal P}^d_{\min (n,a,k,b, m )} }
             \Hom(\Vm, \Vn)
  \end{array}
 \end{eqnarray*}
Furthermore, the map in the left column induced by the convolution
is identified with the map in the right column given by the
obvious composition.
\section{The $\pair$-duality and Schur duality}
      \label{sec_schur}
   Denote by $\B$ the flag manifold of $GL_d$, namely
$$
 {\B} :=
  \{{\frak F} = (0 = F_0 \subset F_1 \subset \ldots \subset F_d = \C^d ) |
       \dim {\cal F}_i / {\cal F}_{i -1} =1, i =1, \ldots, d \}.
$$
Recall that $\cal N$ is the nilpotent cone of $gl_d$. We define
$$
 \widetilde{\cal N} : =
  \{ (x, {\frak F}) \in {\cal N} \times {\B} \mid
     x (F_i) \subset F_{i -1},  i =1, \ldots, d \}.
$$
The following diagram is commutative and $GL_d$-equivariant:
\begin{eqnarray*}
 \begin{array}{rcl}
    &  \widetilde{\cal N} &   \\
  \mu\swarrow& & \searrow \\
  {\cal N } \quad& & \quad\B
 \end{array}
\end{eqnarray*}

We also have a natural $GL_d$-equivariant vector bundle
isomorphism between $\widetilde{\cal N}$ and the
cotangent bundle $T^*{\B}$.
The projection to $\cal N$ is called the Springer resolution.
${\B}_x := \mu^{-1} (x) $, $x \in \cal N$ is the Springer fiber
(cf. e.g. \cite{CG}).

We define the fibred products
\begin{eqnarray*}
    Z  := \widetilde{\cal N} \times_{\cal N} \widetilde{\cal N}, \quad
   \Wn  := \Mn \times_{\Nn} \widetilde{\cal N},
\end{eqnarray*}
and denote
\begin{eqnarray*}
 \Wn \hat{\circ} \Wm :=
   (\Mn \times_{\Nn} \widetilde{\cal N}) \circ
   (\widetilde{\cal N} \times_{\Nn} \Mm).
\end{eqnarray*}
One easily shows that
\begin{eqnarray*}
  \Wn \circ Z = \Wn, \quad \znn \circ \Wn = \Wn,
  \quad \Wn \hat{\circ} \Wm = \znm.
\end{eqnarray*}

The algebra structure on $H(Z)$
is given as a special case of Proposition~\ref{prop_alg} below.
\begin{theorem} \label{th_springer} {\em (Springer)}
   The algebra  $H(Z)$ is isomorphic to the group algebra $\C [S_d] $
  of the symmetric group $S_d$.
\end{theorem}

\begin{theorem}   \label{th_schur}
 {\em (Schur duality)}
   $H(\Wn)$ is isomorphic to the $d$-th tensor
   space of $\C^n$, such that the following diagram commutes:
 \begin{eqnarray*}
  \begin{array}{ccccc}
      U(gl_n) & \looparrowright
              & (\C^n )^{\bigotimes d}   & \looparrowleft & \C [S_d] \\
   \downarrow &
              & \parallel             &  & \parallel \\
      H(\znn) & \looparrowright
              & H(\Wn)                     & \looparrowleft & H(Z)  \\
    \parallel &
              & \parallel             &  & \parallel \\
  \bigoplus_{\lambda \in {\cal P}^d_n }\End (\Vn) & \looparrowright
       & \bigoplus_{\lambda \in {\cal P}^d_{n} }\Vn \bigotimes S_{\lambda}
       & \looparrowleft & \bigoplus_{\lambda \in {\cal P}^d}\End (S_{\lambda})
  \end{array}
 \end{eqnarray*}
 where $S_{\lambda}$ is the irreducible representation of $S_d$
 parameterized by $\lambda$.
\end{theorem}

The variety $\Wn$ was earlier introduced by Lusztig and Grojnowski
\cite{GL}. Although we cannot find the Schur duality in terms of
Borel-Moore homology written down explicitly anywhere, it is
certainly well known to experts (see \cite{GL, GRV}). Next, we
shall derive Theorem~\ref{th_schur} and also Springer's theorem
from the $\pair$-duality. Given a $gl_d$-module $U$, we shall use
$U^{ {\frak h}_d, det}$ to denote the zero-weight space in $U$ of
weight $(1,1, \ldots, 1)$ with respect to  the Cartan subalgebra
${\frak h}_d$.

Set $m =d$. Note that $\B$ is a connected component of ${}^d{\cal
F}$. By Remark~\ref{rem_weight}, We see that the irreducible
components in $\fn \times {}^d {\cal F}$ whose fundamental classes
are of zero weight $(1, \ldots, 1)$ with respect to ${\frak h}_d$
are those coming from $\fn \times_{ {\cal N} } \B$. In this way we
see that
\begin{eqnarray}  \label{eq_imply}
    H(\Wn) = H(\znd)^{ {\frak h}_d, det}.
\end{eqnarray}
In view of the identification (\ref{eq_isom}) and
Remark~\ref{rem_weight}, this is a geometric interpretation of the
following isomorphism of $(gl_n , S_d)$-modules (cf. \cite{H}):
\begin{eqnarray}  \label{eq_equiva}
  (\C^n)^{\bigotimes d} =
  \left(
         S^d (\C^n \bigotimes \C^d )
  \right)^{ {\frak h}_d, det}.
\end{eqnarray}
The $\pair$-duality identifies the r.h.s. of Eqs.~(\ref{eq_imply}) and
(\ref{eq_equiva}), and so it follows that $ H(\Wn)= (\C^n)^{\bigotimes d}.$

Note that ${\B}_{x_{\lambda} } \subset {}^d {\cal F}_{x_{\lambda} }$.
An argument similar to what leads to Eq.~(\ref{eq_imply}) shows that
\begin{eqnarray}  \label{eq_same}
  H({\cal B}_{x_{\lambda}})  =
   H({}^d {\cal F}_{x_{\lambda}})^{ {\frak h}_d, det}
 = {}^d V_{\lambda}^{ {\frak h}_d, det}.
\end{eqnarray}

We note that the Weyl group $S_d$ acts on the zero-weight vector space of
a $gl_d$ module. The following is a well-known
fact, cf. \cite{H} and references therein.

\begin{proposition}   \label{prop_known}
   For $\lambda \in {\cal P}^d_n$, as an $S_d$-module
  ${}^d V_{\lambda}^{ {\frak h}_d, det}$ is
  irreducible and isomorphic to $S_{\lambda}$.
\end{proposition}

By combining Theorem~\ref{th_gldual}, the isomorphism (\ref{eq_imply}) and
Proposition~\ref{prop_known}, we obtain
\begin{eqnarray*}
  H(\Wn) & =& H(\znd)^{ {\frak h}_d, det}  \\
   & =&  \bigoplus_{\lambda \in {\cal P}^d_n }
         H({}^d {\cal F}_{x_{\lambda}}) \bigotimes
          H({\cal B}_{x_{\lambda}})    \\
   & =&  \bigoplus_{\lambda \in {\cal P}^d_n }
          \Vn \bigotimes \Vm^{ {\frak h}, det}   \\
   & =&  \bigoplus_{\lambda \in {\cal P}^d_n }
          \Vn \bigotimes S_{\lambda}.
\end{eqnarray*}

To complete the diagram of Schur duality in Theorem~\ref{th_schur}
it only remains to prove Theorem~\ref{th_springer}.

\begin{demo}{Proof of Theorem~\ref{th_springer}}
  Let us consider the action
 of ${\frak h}_n$ on $H(\zdd)^{ {\frak h}_d, det}$ in the case $n =d$.
 Noting that $Z \subset {}^d Z^d,$ we obtain an analog of (\ref{eq_imply}):
 \begin{eqnarray}   \label{eq_double}
    H(\zdd)^{ {\frak h}_d \bigoplus {\frak h}_d, det \bigoplus det}= H(Z).
 \end{eqnarray}
 By combining the $(gl_d, gl_d)$-duality(Theorem~\ref{th_gldual}),
 the isomorphism (\ref{eq_double})
 and Proposition~\ref{prop_known}, we have
 $$
  H(Z)  \approx \bigoplus_{\lambda \in {\cal P}^d }
  S_{\lambda} \bigotimes S_{\lambda}  \approx \C [S_d ].
 $$
\end{demo}

Indeed the Schur duality also implies the $\pair$-duality.
The equivalence between $\pair$-duality and
Schur duality was established in a completely different setting \cite{H}.
However we hope that the reader may find it illuminating
to see how the interplay between the $\pair$-duality
and the Schur duality is reflected in our geometric setup.

Let us assume the Schur duality. $\Wn$ is a union of
\begin{eqnarray}     \label{eq_bridge}
 \Wn ({\cal O}) :=
 GL_d \times_{ G(x)} ( \fn_x \times {\B}_x)
        \quad (x \in {\cal O})
\end{eqnarray}
over all nilpotent orbits $\cal O$ in $\Nn$.
The composition of set between $ \Wn ({\cal O})$ and
$ \Wm ( {\cal O}^{'})$ (in the sense of $\Wn \hat{\circ} \Wm$)
is non-zero if and only if ${\cal O}^{'} = \cal O$, and
\begin{eqnarray} \label{eq_comp}
  \Wn ({\cal O}) \hat{\circ} \Wm ({\cal O} )
  = GL_d \times_{ G(x)} ( \fn_x \times {\fm}_x)
  = \znm ({\cal O} ).
\end{eqnarray}

The Schur duality states that
\begin{eqnarray*}
 H(\Wn) = \bigoplus_{\lambda \in {\cal P}^d_{n} } H(\fn_{x_{\lambda}})
 \bigotimes H({\B}_{x_{\lambda}})
 = \bigoplus_{\lambda \in {\cal P}^d_{n} }\Vn \bigotimes S_{\lambda}.
\end{eqnarray*}
It follows from (\ref{eq_bridge})
that the isotropic space
$ H(\fn_{x_{\lambda}}) \bigotimes H({\B}_{x_{\lambda}}) $
is spanned by the fundamental classes of irreducible
components of $ \Wn ({\cal O}_{\lambda})$, where
${\cal O}_{\lambda} $ denotes the orbit
of $x_{\lambda}$. In particular $H({\fn}_{x_{\lambda}})$
is isomorphic to $\Vn$. Then Eq.~(\ref{eq_comp}) and
the fact $\znm = \bigsqcup_{ \cal O} \znm ({\cal O})$ implies
the $\pair$-duality
\begin{eqnarray*}
 H (\znm)   = \bigoplus_{\lambda \in {\cal P}^d_{\min (n,m)} }
     H(\fn_{x_{\lambda}}) \bigotimes  H(\fm_{x_{\lambda}})
       = \bigoplus_{\lambda \in {\cal P}^d_{\min (n,m)} }
     \Vn \bigotimes \Vm .
\end{eqnarray*}

The rest of this section is the counterpart of results
in Section~\ref{sec_main} by substituting $\fn, \Nn$ etc with
${\B}, \cal N$ etc. The proofs are totally analogous which we will omit.

We introduce a new variety
$$
 Z_k := \widetilde{\cal N} \times_{\Nk} \widetilde{\cal N} \quad (k \geq 0).
$$
We have the composition of sets
$$
  Z_k \circ Z_l = Z_{ \min (k, l)}.
$$
We can prove as in the preceding section that
$Z_k$ is a lagrangian subvariety of
$\widetilde{\cal N} \times_{\cal N} \widetilde{\cal N}$.
It follows from the dimension property of the convolution
in Borel-Moore homology that the convolution induces a map
$$
  \star : H (Z_k ) \bigotimes H(Z_l ) \longrightarrow H(Z_{ \min (k, l)}).
$$
In particular we obtain an analog of Proposition~\ref{prop_spec}.
Note that $Z_k \subset Z$, and $Z_k = Z$ iff $k \geq d$.

\begin{proposition}  \label{prop_alg}
 $(H(Z_k), \star) $ is an associative algebra with unit. $H(Z_k)$ carries a
 $H(Z)$ bi-module structure.
\end{proposition}

One easily shows that
$$
 Z \circ {\B}_x = {\B}_x, \quad {\B}_x \circ Z = {\B}_x.
$$
The dimension property of $Z$ ensures that the convolution induces
a $H(Z)$ bi-module structure on $H({\B}_x)$. The following
theorem is an analog of Theorem~\ref{th_duality}.
In the case when $k \geq d$ so that $Z_k = Z$ it
reduces to Springer's theorem.

\begin{theorem}  \label{th_easy}
 The algebra $  H(Z_k )$
  is isomorphic to $ \bigoplus_{ \lambda \in {\cal P}^d_k}
     \End ( H({\B}_{x_{\lambda}} ))$.
   This is also an $S_d$-bi-module isomorphism.
\end{theorem}

In view of Eq.~(\ref{eq_same})
and Proposition~\ref{prop_known}, we reformulate
Theorem~\ref{th_easy} as follows.

\begin{theorem}
 The algebra $H(Z_k)$ is isomorphic to
 $\bigoplus_{ \lambda \in {\cal P}^d_k} \End ( S_{\lambda}).$
  This is also an $S_d$-bi-module isomorphism.

\end{theorem}

\noindent{\bf Acknowledgment.} This paper was initiated during my
stay at Max-Planck-Institut f\"ur Mathematik at Bonn, and was
first presented in the Lie group seminar at Yale University in
1998. It is a pleasure to thank MPI for its warm hospitality and
stimulating atmosphere.

Department of Mathematics, North Carolina State University,
Raleigh, NC 27695. E-mail: wqwang@math.ncsu.edu

\end{document}